\documentclass[reqno,12pt,dvips]{article}
\usepackage{amsmath,amsfonts,amssymb,amsthm,enumerate,graphicx,ifthen,bm,geometry,needspace,fullpage}
\usepackage[bookmarks=true,implicit=true]{hyperref}
\usepackage[all]{hypcap}

\newcommand{\smallfigwid}{200pt}
\newcommand{\mediumfigwid}{300pt}
\newcommand{\largefigwid}{400pt}

\renewcommand{\exp}[1]{e^{#1}}
\newcommand{\abs}[1]{\left|#1\right|}
\newcommand{\rb}[1]{\left(#1\right)}
\newcommand{\sqb}[1]{\left[#1\right]}
\newcommand{\setb}[1]{\left\{#1\right\}}
\newcommand{\mymbox}[1]{\quad\mbox{#1}\quad}
\newcommand{\st}{~:~}
\newcommand{\vep}{\varepsilon}
\newcommand{\hrho}{\widehat{\rho}}
\newcommand{\hsigma}{\widehat{\sigma}}
\newcommand{\cM}{\mathcal{M}}
\newcommand{\cO}{\mathcal{O}}
\newcommand{\sfrac}[2]{\textstyle\frac{#1}{#2}\displaystyle}
\newcommand{\deriv}[2]{\frac{d#1}{d#2}}
\newcommand{\derivb}[2]{\frac{d}{d#2}\rb{#1}}
\newcommand{\pderiv}[2]{\frac{\partial#1}{\partial#2}}
\newcommand{\integral}[4]{\int_{#1}^{#2} \, #3 \, d#4}
\newcommand{\normal}{\widehat{\vec{n}}}
\newcommand{\Rd}[1]{\mathbb{R}^{#1}}
\newcommand{\Nd}[1]{\mathbb{N}^{#1}}
\newcommand{\Qd}[1]{\mathbb{Q}^{#1}}
\newcommand{\Zd}[1]{\mathbb{Z}^{#1}}
\renewcommand{\vec}[1]{\textbf{#1}}
\newcommand{\hvec}[1]{\widehat{\vec{#1}}}
\newcommand{\mat}[1]{\textbf{#1}}
\newcommand{\dotprod}{\bullet}
\newcommand{\norm}[1]{\left\|#1\right\|}
\newcommand{\twobyonematrix}[2]{\begin{pmatrix}#1\\#2\end{pmatrix}}
\newcommand{\stwobyonematrix}[2]{\rb{\begin{smallmatrix}#1\\#2\end{smallmatrix}}}
\newcommand{\twobytwomatrix}[4]{\begin{pmatrix}#1&#2\\#3&#4\end{pmatrix}}
\newcommand{\stwobytwomatrix}[4]{\rb{\begin{smallmatrix}#1&#2\\#3&#4\end{smallmatrix}}}
\newcommand{\seq}[4]{\setb{{#1}_{#2}}_{\ifthenelse{\equal{#3}{}}{}{#2=#3}}^{#4}}
\newcommand{\floor}[1]{\left\lfloor#1\right\rfloor}

\theoremstyle{plain}
\newtheorem{thm}{\bf{Theorem}}
\newtheorem{cor}[thm]{\bf{Corollary}}
\newtheorem{lem}[thm]{\bf{Lemma}}
\newtheorem{prop}[thm]{\bf{Proposition}}

\theoremstyle{definition}
\newtheorem{defn}[thm]{\bf{Definition}}
\newtheorem{rem}[thm]{\bf{Remark}}
\newtheorem{rems}[thm]{\bf{Remarks}}

\renewenvironment{proof}{\medskip\noindent\em Proof:
\rm}{\hspace*{\fill}$\square$\medskip}

\begin{document}

\title{Properties of the Michaelis-Menten Mechanism in\\Phase Space\\\text{}\\Matt S. Calder, David Siegel\\\text{}\\\small{Department of Applied Mathematics, University of Waterloo}\\\small{200 University Avenue West,
Waterloo, Ontario N2L 3G1, Canada}}

\date{}

\maketitle

\textbf{Abstract:} We study the two-dimensional reduction of the
Michaelis-Menten reaction of enzyme kinetics.  First, we prove the
existence and uniqueness of a slow manifold between the horizontal
and vertical isoclines.  Second, we determine the concavity of all
solutions in the first quadrant.  Third, we establish the asymptotic
behaviour of all solutions near the origin, which generally is not
given by a Taylor series.  Finally, we determine the asymptotic
behaviour of the slow manifold at infinity.  To this end, we show
that the slow manifold can be constructed as a centre manifold for a
fixed point at infinity.


\section{Introduction}

The Michaelis-Menten (more accurately the Michaelis-Menten-Henri)
mechanism is the simplest chemical network which models the
formation of a product through an enzymatic catalysis of a
substrate.  See, for example,
\cite{FraserRoussel1991,Henri,MichaelisMenten,SegelSlemrod},
Chapter~1 of \cite{KeenerSneyd}, and Chapter~10 of \cite{LinSegel}.
In particular, an enzyme reacts with the substrate and reversibly
forms an intermediate complex, which then decays into the product
and original enzyme. Symbolically,
\[
    S + E \overset{k_1}{\underset{k_{-1}}{\rightleftarrows}} C \overset{k_2}{\rightarrow} P + E,
\]
where $S$ stands for substrate, $E$ stands for enzyme, $C$ stands
for complex, and $P$ stands for product.  By the Law of Mass Action
we have a set of four differential equations for the concentrations
$s$, $e$, $c$, and $p$:
\begin{align*}
    \dot{s} &= k_{-1} c - k_1 s e, \\
    \dot{e} &= ( k_{-1} + k_2 ) c - k_1 s e, \\
    \dot{c} &= k_1 s e - ( k_{-1} + k_2 ) c, \\
    \dot{p} &= k_2 c.
\end{align*}
Taking initial conditions ${s(0)=s_0}$, ${e(0)=e_0}$, ${c(0)=0}$,
and ${p(0)=0}$, we have ${e=e_0-c}$ and two independent differential
equations:
\begin{align*}
    \dot{s} &= k_{-1} c - k_1 s ( e_0 - c ), \\
    \dot{c} &= k_1 s ( e_0 - c ) - ( k_{-1} + k_2 ) c.
\end{align*}

To simplify even further, the Quasi-Steady-State Approximation
(QSSA) was introduced by Briggs and Haldane\cite{BriggsHaldane}. This takes ${\dot{c}=0}$ to hold after a
short time, giving
\[
    c = \frac{ e_0 s }{ K_m + s }
    \mymbox{and}
    - \! \dot{s} = \dot{p} = \frac{ k_2 e_0 s }{ K_m + s },
\]
where ${K_m := \frac{ k_{-1} + k_2 }{ k_1 }}$ is called the
Michaelis constant. This is generally thought to be valid when ${e_0
\ll s_0}$. The expression for $\dot{p}$ gives a measure of the
velocity of the reaction.  Experiments have been fitted to the
quasi-steady-state approximation. This type of approximation is
often used to simplify other chemical kinetics systems including
those of different and more complicated enzyme reactions which may
involve inhibition or cooperativity
effects\cite{ConradTyson,KeenerSneyd,KlippHerwigKowaldWierlingLehrach}.

Another approximation, though less common than QSSA, is the (rapid)
Equilibrium Approximation (EA), which originated with Henri and was
popularized by Michaelis and Menten.  This takes ${\dot{s}=0}$ to
hold after a short time, giving
\[
    c = \frac{e_0 s}{K_s + s}
    \mymbox{and}
    \dot{p} = \frac{k_2 e_0 s}{K_s + s},
\]
where ${K_s := \frac{k_{-1}}{k_1}}$.

The validity of QSSA can be examined from several points of view.
First the equations are written in a non-dimensional form,
introducing a small positive parameter $\vep$.  This can be done is
several different ways\cite{delaSelvaPinaGarcia-Colin,SegelSlemrod}.
In the traditional scaling\cite{LinSegel,Murray,OMalley},
${\vep=\frac{e_0}{s_0}}$. Following
\cite{Roussel1994,Roussel1997,RousselFraser2001}, we have
\begin{equation} \label{eq001}
    \dot{x} = -x + ( 1 - \eta ) y + x y,
    \quad
    \dot{y} = \vep^{-1} ( x - y - x y ),
\end{equation}
where ${x:=\sfrac{k_1 s}{k_{-1}+k_2}}$ is a scaled substrate
concentration, ${y:=\sfrac{c}{e_0}}$ is a scaled complex
concentration, ${t := k_1 e_0 \tau}$ is a scaled time, $\tau$ is the
original time, $\dot{}=\deriv{}{t}$, and ${e_0:=e(0)+c(0)}$. The
parameters are given by ${\vep:=\sfrac{k_1 e_0}{k_{-1}+k_2}}$ and
${\eta:=\sfrac{k_2}{k_{-1}+k_2}}$. Note that ${\vep
> 0}$ and ${0<\eta<1}$.

When $\vep$ is small, the system \eqref{eq001} is a singular
perturbation problem.  A matched asymptotic analysis yields the QSSA
as the zeroth-order term in the outer expansion\cite{LinSegel}.  The
correctness of this analysis is proved using Tikhonov-Levinson
theory\cite{OMalley}. Explicit bounds on the approximations have
been obtained for small $\vep$\cite{SegelSlemrod}. Centre manifold
theory\cite{Carr} and geometric singular perturbation
theory\cite{Kaper} have been applied to give an invariant manifold
$\cM_\vep$, called a slow manifold, within distance $\cO(\vep)$ of
$\cM_0$, the quasi-steady-state manifold ${y=\frac{x}{1+x}}$.
Trajectories approach the slow manifold exponentially fast and then
evolve along it at a slower rate.

Several chemists have observed and theoretically investigated slow
manifolds which attract other solutions.  In general slow manifolds
are not uniquely defined.  In two-dimensional cases, Fraser and
Roussel\cite{Fraser1988,Roussel1994,Roussel1997,RousselFraser2001}
take as a slow manifold the solution between the horizontal and
vertical isoclines, which are the quasi-steady-state and the
equilibrium approximations, respectively. Roussel, for example,
provided a heuristic argument based on antifunnel theory that there
is indeed such a solution\cite{Roussel1997}. Davis and
Skodje\cite{DavisSkodje} take as a slow manifold the trajectory
joining a saddle at infinity and a stable node, approaching in a
slow direction. This paper was motivated by the work of Fraser and
Roussel.

Occasionally, we may refer to the system \eqref{eq001} in the
compact form
\begin{equation} \label{eq002}
    \dot{\vec{x}} = \vec{g}(\vec{x}).
\end{equation}
We will also work with the one-dimensional version of \eqref{eq001},
given by
\begin{equation} \label{eq003}
    y' = f(x,y),
\end{equation}
where $'=\deriv{}{x}$. Explicitly,
\[
    \vec{g}(\vec{x}) := \twobyonematrix{ -x + ( 1 - \eta ) y + x y }{ \vep^{-1} ( x - y - x y ) }
    \mymbox{and}
    f(x,y) := \frac{ x - y - x y }{ \vep [ -x + ( 1 - \eta ) y + x y ] }.
\]

In this paper, we do not need to assume that $\vep$ is small. The
focus is on the behaviour of solutions in the phase plane, that is,
considering $y$ as a function of $x$. In \S\ref{sec002}, we give the
basic phase portrait of \eqref{eq001} in the first quadrant and the
linearization at the origin.  In \S\ref{sec003}, we describe the
isocline structure which is exploited in subsequent sections. In
\S\ref{sec004}, we prove the existence and uniqueness of the slow
manifold, which we denote by $\cM$, between the horizontal and
vertical isoclines.  These were discussed in a more informal way by
Fraser (see, for example, \cite{Fraser1988}). In \S\ref{sec005}, we
determine the concavity of all solutions except the slow manifold by
analyzing an auxiliary function.  In \S\ref{sec006}, we determine
the behaviour of solutions near the origin by using Poincar\'{e}'s
Theorem (see, for example, \cite{Arnold} p.190); one-dimensional
solutions $y(x)$ are generally not given by a Taylor series, which
has sometimes been assumed. This analysis applies to any
two-dimensional system with a Hurwitz-stable equilibrium point. In
\S\ref{sec007}, we determine when solutions enter $\Gamma_1$, which
is a region bounded below by the horizontal isocline and above by
the isocline for the slope of the slow manifold at the origin. In
\S\ref{sec008}, we establish properties of the slow manifold:
concavity, monotonicity, and asymptotic behaviour at the origin and
infinity. Finally, in \S\ref{sec009} we state some open questions.

\section{Phase Portrait} \label{sec002}

\begin{figure}[t]
\begin{center}
    \includegraphics[width=\largefigwid]{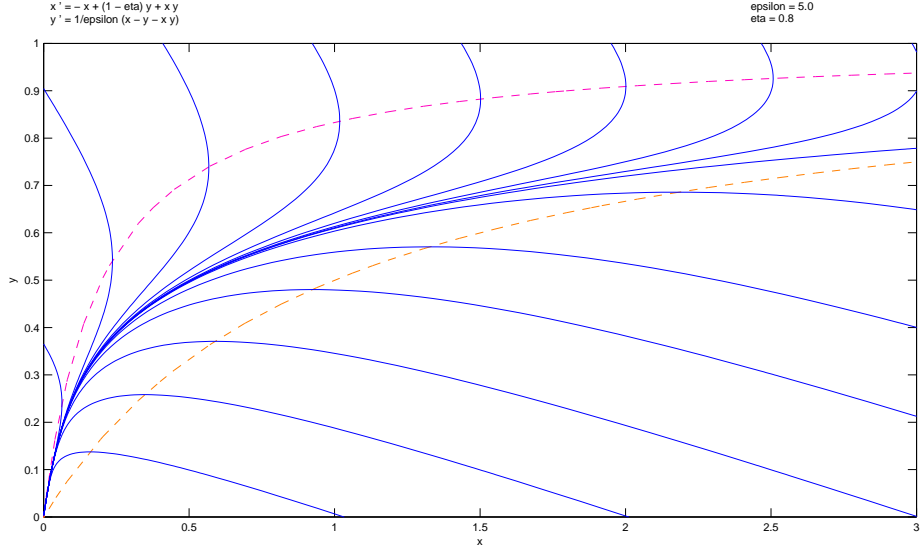}
    \caption{A phase portrait for \eqref{eq001} for ${\vep=5.0}$
    and ${\eta=0.8}$.} \label{fig001}
\end{center}
\end{figure}

The qualitative behaviour of solutions is revealed by the phase
portrait. See, for example, Figure~\ref{fig001}, which is a phase
portrait for certain values of the parameters. To find the
horizontal and vertical isoclines, set, respectively, ${\dot{y}=0}$
and ${\dot{x}=0}$ in \eqref{eq001} to obtain the graphs
\[
    y = H(x) := \frac{ x }{ 1 + x }
    \mymbox{and}
    y = V(x) := \frac{ x }{ 1 - \eta + x }.
\]
Note that the EA corresponds to the vertical isocline and the QSSA
corresponds to the horizontal isocline. (One may refer to ${y=H(x)}$
as the quasi-steady-state manifold and ${y=V(x)}$ as the rapid
equilibrium manifold.) Observe that
\[
    \lim_{x \to \infty} H(x) = 1 = \lim_{x \to \infty} V(x)
    \mymbox{and}
    H(0) = 0 = V(0).
\]
Since both $H$ and $V$ are strictly increasing and ${V(x)>H(x)}$ for
all ${x>0}$, there is a narrow region between the isoclines:
\[
    \Gamma_0 := \setb{ (x,y) \st x > 0, ~ H(x) \leq y \leq V(x) }.
\]

\begin{thm}
Let $\vec{x}(t)$ be a solution to \eqref{eq001}.  Then, there exists
a ${t^*>0}$ such that
\[
    (x(t),y(t)) \in \Gamma_0 ~ \forall ~ t \geq t^*.
\]
\end{thm}

\begin{proof}
First, we show that solutions can enter $\Gamma_0$ but not leave it,
that is, show that $\Gamma_0$ is positively invariant.  It follows
from the differential equation that ${\vec{g} \dotprod \normal < 0}$
along both the vertical and horizontal isoclines, where $\normal$ is
the outward unit normal vector.  Hence, $\Gamma_0$ is positively
invariant.

Second, we establish that solutions outside $\Gamma_0$ eventually
enter $\Gamma_0$.  Call $y(x)$ the corresponding one-dimensional
solution. If $y(x)$ is below the horizontal isocline $H(x)$, then
${-\vep^{-1} \leq y'(x) < 0}$, and so $y$ must intersect $H$ for a
lower value of $x$. Similarly, if $y(x)$ is above the vertical
isocline $V(x)$, then ${-\infty<y'(x)<-\vep^{-1}}$, and so $y$ must
intersect $V$ for a higher value of $x$.
\end{proof}

Behaviour of solutions near the origin, the only equilibrium point,
is governed by the linearization matrix
\begin{equation} \label{eq004}
    \mat{A}
    := \left.\pderiv{\vec{g}}{\vec{x}}\right|_{ \vec{x} = \vec{0} }
    = \twobytwomatrix{-1}{1-\eta}{\vep^{-1}}{-\vep^{-1}}.
\end{equation}
The eigenvalues of $\mat{A}$ are given by
\begin{equation} \label{eq005}
    \lambda_\pm := \frac{ -( \vep + 1 ) \pm \sqrt{ ( \vep + 1 )^2 - 4 \vep \eta } }{ 2 \vep },
\end{equation}
which are real-valued and distinct. One can prove that
\[
    \lambda_- < -1 < \lambda_+ < 0,
\]
thus implying that the origin is asymptotically stable.
Corresponding eigenvectors are
\[
    \vec{v}_\pm
    :=
    \stwobyonematrix{ 1 - \eta }{ \lambda_\pm + 1 }.
\]
Observe that $\vec{v}_+$ points into the positive quadrant while
$\vec{v}_-$ does not.  The slope of the eigenvector $\vec{v}_+$ at
the origin is very important, and will be denoted by
\[
    \sigma := \frac{ \lambda_+ + 1 }{ 1 - \eta }.
\]
Asymptotically, we have that
\[
    \sigma = 1 + \vep \eta + \cO\rb{\vep^2}
    \mymbox{as}
    \vep \to 0.
\]
Observe also that the slope of the slow manifold at the origin lies
between the slope of the horizontal isocline and the slope of the
vertical isocline. That is,
\[
    1 < \sigma < ( 1 - \eta )^{-1}.
\]

The original, time-dependent differential equation \eqref{eq001} has
linearization
\begin{equation} \label{eq006}
    \dot{\vec{x}} = \mat{A} \vec{x},
\end{equation}
where the matrix $\mat{A}$ is as in \eqref{eq004}. Since the
eigenvalues are real-valued and distinct, the initial value problem
\[
    \dot{\vec{x}} = \mat{A} \vec{x},
    \quad
    \vec{x}(0) = \vec{x}_0
\]
has solution of the form
\[
    \vec{x}(t) = c_- \exp{\lambda_- t} \vec{v}_- + c_+ \exp{\lambda_+ t} \vec{v}_+.
\]
We will assume, to avoid triviality, that ${\vec{x}_0 \ne \vec{0}}$.
The coefficients $c_\pm$ can be determined in terms of left
eigenvalues and eigenvectors.  The left eigenvectors can be taken to
be
\[
    \hvec{v}_\pm := \stwobyonematrix{ \vep^{-1} }{ \lambda_\pm + 1 }.
\]
Using the orthogonality condition ${\hvec{v}_\pm^T \vec{v}_\mp = 0}$
we see that
\[
    c_\pm = \frac{ \hvec{v}_\pm^T \vec{x}_0 }{ \hvec{v}_\pm^T \vec{v}_\pm }.
\]

\begin{prop}
Let $y$ be a solution to \eqref{eq003} which lies inside $\Gamma_0$
for ${x\in(0,a)}$, where ${a>0}$. Then,
\[
    \lim_{x \to 0^+} y(x) = 0
    \mymbox{and}
    \lim_{x \to 0^+} y'(x) = \sigma.
\]
\end{prop}

\begin{proof}
We should begin by emphasizing that solutions $\vec{x}(t)$ to
\eqref{eq001} enter and forever remain in the interior of
$\Gamma_0$. By hypothesis, ${H(x)<y(x)<V(x)}$ for all ${x\in(0,a)}$.
The Squeeze Theorem establishes the first limit since ${H(0)=0}$ and
${V(0)=0}$.

To establish the second result, observe that the function $\vec{g}$,
as in \eqref{eq002}, is of class $C^2$ with
${\vec{g}(\vec{0})=\vec{0}}$ and the matrix $\mat{A}$ has strictly
negative eigenvalues.  It follows from Hartman's Theorem (see, for
example, \cite{Perko} p.127), which is a stronger version of the
Hartman-Grobman Theorem and applies even in cases of resonance, that
the phase portrait of \eqref{eq001} behaves like the phase portrait
of \eqref{eq006} diffeomorphically in a neighbourhood of the origin.
Therefore, solutions to the nonlinear system have slope $\sigma$ as
they approach the origin too.
\end{proof}

\begin{rem}
Since any solution $\vec{x}(t)$, except the trivial solution,
eventually enters $\Gamma_0$ and then approaches the origin
asymptotically, we can now say definitively that the origin is
globally asymptotically stable.
\end{rem}

\section{The Isocline Structure} \label{sec003}

The nature of the level curves, or isoclines, ${c=f(x,y(x))}$
reveals to us a surprising amount of insight into the behaviour of
solutions to \eqref{eq003}.  Consider
\begin{equation} \label{eq007}
    c = f(x,y(x)),
\end{equation}
where $f$ is as in the differential equation \eqref{eq003} and
${c\in\Rd{}}$.  For a given ${x>0}$ and ${c\in\Rd{}}$, \eqref{eq007}
is invertible and we can solve for $y(x)$, yielding
\begin{equation} \label{eq008}
    y(x) = F(x,c) := \frac{ x }{ K(c) + x },
\end{equation}
for ${c \ne -\vep^{-1}}$, where
\[
    K(c) := \frac{ 1 + \vep ( 1 - \eta ) c }{ 1 + \vep c }.
\]
Observe that ${y'(x) = f(x,y(x))}$ and ${y(x) = F(x,y'(x))}$ for
solutions $y$ of \eqref{eq003} for values of $x$ for which the
solution is defined. Throughout this paper, level curves of $f$ will
be denoted by $w$. If the slope associated with $w$ is required, we
specify this and write ${w(x)=F(x,c)}$. For completeness, we will
agree that ${F\rb{x,-\vep^{-1}} = 0}$.

\begin{figure}[t]
\begin{center}
    \includegraphics[width=\mediumfigwid]{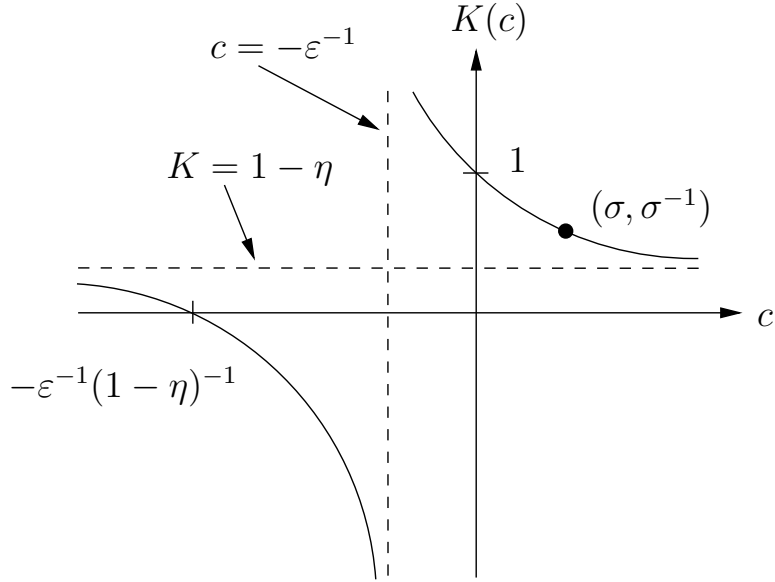}
    \caption{Graph of the function $K(c)$ for arbitrary ${\vep>0}$ and ${\eta\in(0,1)}$.} \label{fig002}
\end{center}
\end{figure}

\begin{figure}[t]
\begin{center}
    \includegraphics[width=\largefigwid]{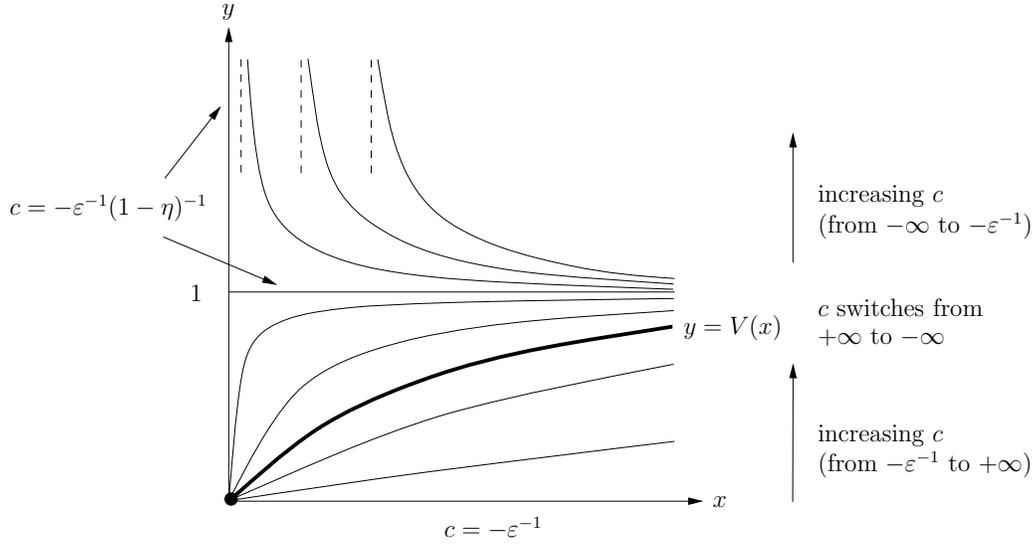}
    \caption{Sketch of the isocline structure of \eqref{eq003} for arbitrary ${\vep>0}$ and ${\eta\in(0,1)}$.} \label{fig003}
\end{center}
\end{figure}

The function $K$, which is sketched in Figure~\ref{fig002}, and the
isoclines, which are sketched in Figure~\ref{fig003}, have the
following, easy-to-prove properties.

\begin{figure}[t]
\begin{center}
    \includegraphics[width=\smallfigwid]{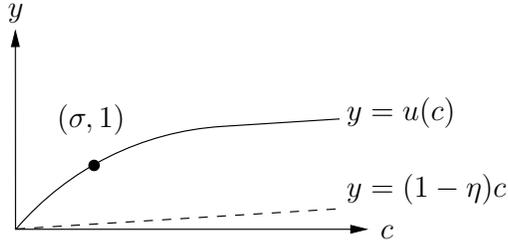}
    \caption{Graph of the function $u(c)$ given in Proposition~\ref{prop001}(d) for arbitrary ${\vep>0}$ and ${\eta\in(0,1)}$.
    The function is strictly increasing for ${c>0}$ for every admissible $\vep$ and $\eta$.  Furthermore, ${u'(0)=1}$ and
    asymptotically $u$ has slope ${1-\eta}$ as ${c\to\infty}$. Incidently, ${u(c)-(1-\eta)c = \vep^{-1}\eta + \cO\rb{c^{-1}}}$ as ${c\to\infty}$.} \label{fig004}
\end{center}
\end{figure}

\begin{prop} \label{prop001}
The isoclines and the function $K$ satisfy the following.
\begin{enumerate}[(a)]
    \item
        The function $K$ is strictly decreasing everywhere except at
        ${c=-\vep^{-1}}$, where it has a vertical asymptote.
    \item
        The function $K$ satisfies ${1 - \eta < K(c) < 1}$ for ${0 < c <
        \infty}$, which corresponds to the interior of $\Gamma_0$. Furthermore,
        we have that ${K(0)=1}$ corresponds to
        the horizontal isocline and ${\lim_{c\to\infty}K(c)=1-\eta}$
        corresponds to the vertical isocline.
    \item
        The function $K$ satisfies the important relation ${K(\sigma) =
        \sigma^{-1}}$.
    \item
        Define the function ${u(c) := c K(c)}$ for ${c>0}$.  The
        function $u$ is strictly increasing, satisfies
        ${u(\sigma)=1}$, and
        \[
            u(c) = c K(c) \in
            \begin{cases}
                (0,1), \mymbox{for} c \in (0,\sigma) \\
                (1,\infty), \mymbox{for} c \in (\sigma,\infty)
            \end{cases}.
        \]
        See Figure~\ref{fig004}.
    \item
        Any of the isoclines $w$ satisfy the
        differential equation
        \begin{equation} \label{eq009}
            w ( w - 1 ) + x w' = 0.
        \end{equation}
\end{enumerate}
\end{prop}

Note that the isoclines are hyperbolas. There are two exceptional
isoclines, namely ${w(x)=F\rb{x,-\vep^{-1}} = 0}$ and ${
w(x)=F\rb{x,-\vep^{-1}(1-\eta)^{-1}} = 1}$. The vertical isocline,
$V$, also is somewhat of an exceptional case. Approaching it from
below, one encounters increasing $c$ up to $+\infty$.  After passing
through $V$, the slopes increase from $-\infty$.

\section{Existence and Uniqueness of the Slow Manifold}
\label{sec004}

We provide a brief review of fences and antifunnels, which form the
backbone of the existence-uniqueness proof that follows.

Phase spaces of differential equations often exhibit curious curves
and regions known as fences, funnels, and antifunnels.  The best
source of information on funnels and antifunnels is
\cite{HubbardWest}, Chapters~1 and 4.

\begin{defn} \label{defn001}
Let ${I=[a,b)}$ be an interval (where ${a<b\leq\infty}$) and
consider the first-order differential equation $y'=f(x,y)$ over $I$.
Let $\alpha$ and $\beta$ be continuously-differentiable functions
satisfying
\begin{equation} \label{eq010}
    \alpha'(x) \leq f(x,\alpha(x))
    \mymbox{and}
    f(x,\beta(x)) \leq \beta'(x)
\end{equation}
for all ${x \in I}$.
\begin{enumerate}[(a)]
    \item
        The curves $\alpha$ and $\beta$ satisfying \eqref{eq010} are, respectively, a \emph{lower fence} and an \emph{upper fence}. If
        there is always a strict inequality in \eqref{eq010}, the fences are
        \emph{strong}.  Otherwise, the fences are \emph{weak}.
    \item
        If ${\beta(x)<\alpha(x)}$ on $I$, then the set
        \[
            \Gamma := \setb{ (x,y) \st x \in I, ~ \beta(x) \leq y \leq \alpha(x) }
        \]
        is called an \emph{antifunnel}.  The antifunnel is \emph{narrowing} if
        \[
            \lim_{x \to b^-} \abs{ \alpha(x) - \beta(x) } = 0.
        \]
\end{enumerate}
\end{defn}

\begin{thm}[Antifunnel Theorem, \cite{HubbardWest} p.31-33] \label{thm001}
Let $\Gamma$ be an antifunnel with strong lower and upper fences
$\alpha$ and $\beta$, respectively, for the differential equation
${y'=f(x,y)}$ over the interval ${I:=[a,b)}$ (where
${a<b\leq\infty}$). Then, there exists a solution $y(x)$ to the
differential equation such that
\[
    \beta(x) < y(x) < \alpha(x) \mymbox{for all} x \in I.
\]
If, in addition, $\Gamma$ is narrowing and ${\pderiv{f}{y}(x,y) \geq
0}$ in $\Gamma$, then the solution $y(x)$ is unique.
\end{thm}

We cannot use $\Gamma_0$ as an antifunnel (in the sense of
Definition~\ref{defn001}).  The key to our proof is considering the
isocline for slope $\sigma$, the slope of the slow manifold at the
origin. We will call this isocline $\alpha$.
Proposition~\ref{prop001}(c) and \eqref{eq008} tell us that $\alpha$
is given by the simple expression
\begin{equation} \label{eq011}
    \alpha(x) = \frac{ x }{ \sigma^{-1} + x }.
\end{equation}
The defining feature of $\alpha$ is ${\sigma \equiv
f(x,\alpha(x))}$. Moreover, this function $\alpha$ has the
remarkable property that ${\alpha'(0)=\sigma}$. That is, the
isocline for slope $\sigma$ has slope $\sigma$ at the origin. Define
the region
\[
    \Gamma_1 := \setb{ (x,y) \st x > 0, ~ H(x) \leq y \leq \alpha(x) },
\]
which is a subset of $\Gamma_0$ because ${H(x) < \alpha(x) < V(x)}$
for all ${x>0}$.

\begin{thm}\text{} \label{thm002}
\begin{enumerate}[(a)]
    \item
        There exists a unique slow manifold $y=\cM(x)$ in $\Gamma_1$ for the
        differential equation \eqref{eq003}.
    \item
        The slow manifold ${y=\cM(x)}$ is also the only solution
        that lies entirely inside $\Gamma_0$.
\end{enumerate}
\end{thm}

\begin{proof}
\begin{enumerate}[(a)]
    \item
        We will show that the Antifunnel Theorem can be applied to
        the interval $[a,\infty)$, where ${a>0}$ is arbitrary. First, we show
        that the curve ${y=\alpha(x)}$ is a strong lower fence
        and the curve ${y=H(x)}$ is a strong upper fence for the differential
        equation \eqref{eq003} for ${x>0}$.  Now, the derivative of solutions
        along the concave-down curve ${y=\alpha(x)}$ is identically $\sigma$.
        Thus,
        \[
            \alpha'(x)
            < \alpha'(0)
            = \sigma
            = f( x, \alpha(x) )
            \quad
            \forall ~ x > 0.
        \]
        Hence, by definition, ${y=\alpha(x)}$ is a strong lower fence for
        ${x>0}$. To show that ${y=H(x)}$ is a strong upper fence for $x>0$,
        consider that
        \[
            f( x, H(x) ) = 0 < H'(x)
            \quad
            \forall ~ x > 0.
        \]

        Second, observe that the strong fences satisfy ${\alpha(x)>H(x)}$ for
        ${x>0}$ and
        \[
            \lim_{x\to\infty} \abs{ \alpha(x) - H(x) } = 0.
        \]
        By definition, $\Gamma_1$ is a narrowing antifunnel.

        Finally, a quick calculation shows that ${\pderiv{f}{y} \geq 0}$ in
        $\Gamma_1$. So, all the conditions for the Antifunnel Theorem
        (Theorem~\ref{thm001}) have been established. Therefore, there
        exists a unique solution ${y=\cM(x)}$ to \eqref{eq003} that lies
        entirely in $\Gamma_1$.
    \item
        Obviously, any solution other than the slow manifold eventually leaves $\Gamma_1$.
        If the solution leaves $\Gamma_1$ through the
        horizontal isocline it also leaves $\Gamma_0$, since both regions share
        the same lower boundary.  If the solution leaves $\Gamma_1$
        through the $\alpha$ isocline, while in $\Gamma_0$ the solution will have slopes
        in the range ${\sigma<y'<\infty}$ and hence will eventually
        leave $\Gamma_0$ since the upper boundary of $\Gamma_0$ is
        bounded above by the line ${y=1}$.
\end{enumerate}
\end{proof}

\needspace{2.0cm}
\begin{rems}~
\begin{enumerate}[(i)]
    \item
        Theorem~\ref{thm002} shows that
        \[
            \frac{x}{1+x} < \cM(x) < \frac{ x }{ \sigma^{-1} + x }
            \quad
            \forall ~ x > 0.
        \]
        Thus, the necessity of the EA is diminished in the sense that $\alpha$ serves as a
        smaller upper bound on $\cM$. Furthermore, it
        follows from the isocline structure that $\cM(x)$ is strictly
        increasing, since solutions of the differential equation inside
        the antifunnel but not on the boundary have strictly positive
        slope.  Note that this bound is especially tight when $\vep$
        is small, since
        \[
            \sigma^{-1} = 1 - \vep \eta + \cO\rb{\vep^2}
            \mymbox{as}
            \vep \to 0.
        \]
    \item
        Slow manifolds, like centre manifolds, are generally not
        unique and are defined locally.  In our case, all solutions that have slope $\sigma$ at the origin
        are slow manifolds. However, we look at the
        global phase portrait and refer to the unique solution
        within $\Gamma_1$ as \emph{the} slow manifold.
\end{enumerate}
\end{rems}

\section{Concavity} \label{sec005}

Let $y$ be a solution to \eqref{eq003}, which we assume is not the
slow manifold because we will deal with that case later. Then, of
course, ${y'(x)=f(x,y(x))}$ and so by the Chain Rule,
\begin{equation} \label{eq012}
    y''(x) = p(x,y(x)) h(x,y(x)),
\end{equation}
where
\[
    p(x,y) := \vep^{-1} \eta \sqb{ -x + ( 1 - \eta + x ) y }^{-2}
\]
and
\begin{equation} \label{eq013}
    h(x,y) := y ( y - 1 ) + x f(x,y).
\end{equation}
The function $p(x,y)$ is positive everywhere except along the
vertical isocline, where it is undefined.  The function
${h(x):=h(x,y(x))}$, the sign of which determines that of $y''(x)$,
has derivative
\begin{equation} \label{eq014}
    h'(x) = 2 y(x) y'(x) + x p(x) h(x),
\end{equation}
where ${p(x):=p(x,y(x))}$. The concavity of all solutions in all
regions of the non-negative quadrant can be deduced using this
auxiliary function $h$. Table~\ref{tab001} summarizes what we will
develop in this section. They are all suggested by the phase
portrait in Figure~\ref{fig001}.
\begin{table}[t]
\begin{center}
\begin{tabular}{|c|c|}
    \hline
    Region & Concavity of Solutions \\ \hline\hline
    $ 0 \leq y < \cM $ & concave down \\ \hline
    $ \cM < y < \alpha $ & concave down, then inflection point, then concave up \\ \hline
    $ \alpha \leq y < V $ & concave up \\ \hline
    $ V < y < 1 $ & concave down \\ \hline
    $ y \geq 1 $ & concave up, then inflection point, then concave down \\ \hline
\end{tabular}
\caption{A summary of the concavity of solutions of \eqref{eq003} in
the non-negative quadrant.} \label{tab001}
\end{center}
\end{table}

\needspace{2.0cm}
\begin{rems}~
\begin{enumerate}[(i)]
    \item
        Let $y$ be a solution to \eqref{eq003} and fix ${x_0>0}$.
        Define ${w(x):=F(x,y'(x_0))}$ to be the isocline through
        $(x_0,y(x_0))$.  By virtue of the isocline structure,
        ${y''(x_0)>0}$ if and only if ${y'(x_0)>w'(x_0)}$ and ${y''(x_0)<0}$ if and only if ${y'(x_0)<w'(x_0)}$.
        Indeed, from \eqref{eq009} and \eqref{eq013},
        \begin{equation} \label{eq015}
            h(x_0) = x_0 [ y'(x_0) - w'(x_0) ],
        \end{equation}
        which confirms this fact.  The similarity of the form of $h(x)$
        and the differential equation \eqref{eq009} that the isoclines
        satisfy is not a coincidence.
    \item
        The function $h$ cannot tell us anything about the
        concavity of solutions at ${x=0}$, not even by taking a
        limit.
\end{enumerate}
\end{rems}

Many of the following proofs will involve the following elementary
lemma, so we single it out here.  We omit the proof in the interest
of space.

\begin{lem} \label{lem001}
Let $I$ be one of the intervals $[a,b]$, $(a,b)$, $[a,b)$, and
$(a,b]$. Suppose that ${\phi \in C(I)}$ is a function having at
least one zero in $I$.
\begin{enumerate}[(a)]
\item
    If ${I=(a,b]}$ or ${I=[a,b]}$, then the function $\phi$ has
    a right-most zero in $I$.  Likewise, if ${I=[a,b)}$ or
    ${I=[a,b]}$, then the function $\phi$ has a left-most zero in $I$.
\item
    If ${\phi \in C^1(I)}$ and ${\phi'(x)>0}$ for every zero
    of $\phi$ in $I$, then $\phi$ has exactly one zero in $I$.
\end{enumerate}
\end{lem}

\begin{prop} \label{prop002}
Let $y(x)$ be any solution to \eqref{eq003} lying below the slow
manifold, say with domain $(0,a]$ and ${y(a)=0}$.  Then, $y$ is
concave down for all ${x \in (0,a]}$.
\end{prop}

\begin{proof}
Let $h$ be defined as in \eqref{eq013} with respect to the solution
$y$. There are two regions to consider, namely where ${y(x)>H(x)}$
and where ${y(x) \leq H(x)}$.  It is clear from \eqref{eq012} and
\eqref{eq013} that ${y''(x)<0}$ for ${y(x) \leq H(x)}$, noting that
${0 \leq y(x) < 1}$ and ${y'(x)<0}$.  The solution $y(x)$ crosses
the horizontal isocline, say at ${x = x_2 \in (0,a)}$. Here,
${h(x_2)<0}$ using \eqref{eq013}.  Suppose that the proposition is
false and that there are one or more inflection points.  Applying
Lemma~\ref{lem001}, let ${x_1 \in (0,x_2)}$ be the right-most zero
of $h$.  Now, from \eqref{eq014},
\[
    h'(x_1) = 2 y(x_1) y'(x_1) > 0.
\]
Then, $h$ is positive in a neighbourhood to the right of $x_1$.
Since ${h(x_2)<0}$, by the Intermediate Value Theorem, $h$ has a
zero in $(x_1,x_2)$ which contradicts the fact that $x_1$ is the
right-most zero. Therefore, there is no inflection point.
\end{proof}

\begin{prop}
Let $y$ be a solution to \eqref{eq003} between $\alpha$ and $V$ over
$(a,b)$.  Then, $y$ is concave up on $(a,b)$.
\end{prop}

\begin{proof}
Fix ${x_0 \in (a,b)}$ and let ${c:=y'(x_0)}$ and ${r:=K(c)}$.  Let
${w(x):=F(x,c)}$ be the isocline through $(x_0,y(x_0))$. With $h$
defined as in \eqref{eq013} with respect to $y$, we have
\[
    h(x_0)
    = x_0 [ y'(x_0) - w'(x_0) ]
    = x_0 \sqb{ \frac{ c ( r + x_0 )^2 - r }{ ( r + x_0 )^2 } },
\]
where we used the expression for $h$ in \eqref{eq015}. Since
${c>\sigma}$, applying Proposition~\ref{prop001} we know ${r c
> 1}$ which implies ${r > \sqrt{ r c^{-1} }}$. Suppose, on the
contrary, that ${y''(x_0) \leq 0}$.  Then,
\[
    c ( r + x_0 )^2 - r \leq 0
    \implies
    x_0 \leq \sqrt{rc^{-1}} - r < 0.
\]
This is a contradiction.
\end{proof}

\begin{prop}
Let $y$ be a solution to \eqref{eq003} lying between $\cM$ and
$\alpha$, with domain $(0,a]$ and ${y(a)=\alpha(a)}$.  Then, $y$ has
exactly one inflection point ${x_1 \in (0,a)}$.  Moreover, $y$ is
concave down on $(0,x_1)$ and concave up on $(x_1,a)$.
\end{prop}

\begin{proof}
We know ${y'(0)=\sigma}$ and ${y'(a)=\sigma}$.  Hence, by Rolle's
Theorem, $y$ has an inflection point ${x_1\in(0,a)}$.  To prove
uniqueness of the inflection point, let $h$ be as in \eqref{eq013}
with respect to the solution $y$.  Now, if $x$ is a zero of $h$,
then
\[
    h'(x) = 2 y(x) y'(x) > 0.
\]
By Lemma~\ref{lem001}, there is at most one zero of $h$. Moreover,
since ${h'(x_1)>0}$, $y$ is concave down on $(0,x_1)$ and concave up
on $(x_1,a)$.

\end{proof}

\begin{prop}
Let $y$ be a solution to \eqref{eq003} which lies above $V$ and
below $1$, with domain $[a,b)$ and ${\lim_{x \to b^-} y(x) = V(b)}$.
Then, $y$ is concave down for all ${x \in [a,b)}$.
\end{prop}

\begin{proof}
It is clear from the expression for $h$, \eqref{eq013}, where $h$ is
defined with respect to the solution $y$, and the fact that ${y'<0}$
in that region.
\end{proof}

\begin{prop}
Let $y$ be a solution to \eqref{eq003} which lies above $1$, with
domain $[0,a]$, where ${y(a)=1}$.  Then, there exists a unique
inflection point ${x_1 \in (0,a)}$.  Moreover, $y$ is concave up
over $[0,x_1)$ and concave down over $(x_1,a]$.
\end{prop}

\begin{proof}
Let $h$ be defined as in \eqref{eq013} with respect to the solution
$y$.  Now,
\[
    y'(0) = -\vep^{-1} ( 1 - \eta )^{-1} = y'(a).
\]
By Rolle's Theorem, there exists ${x_1 \in (0,a)}$ such that
${y''(x_1)=0}$.  The uniqueness of the inflection point follows from
the fact that any zero $x$ of $h$ satisfies ${h'(x)<0}$ and an
application of Lemma~\ref{lem001}. Moreover, since ${h'(x_1)<0}$,
$y$ is concave up on $[0,x_1)$ and concave down on $(x_1,a]$.
\end{proof}

\begin{figure}[t]
\begin{center}
    \includegraphics[width=\smallfigwid]{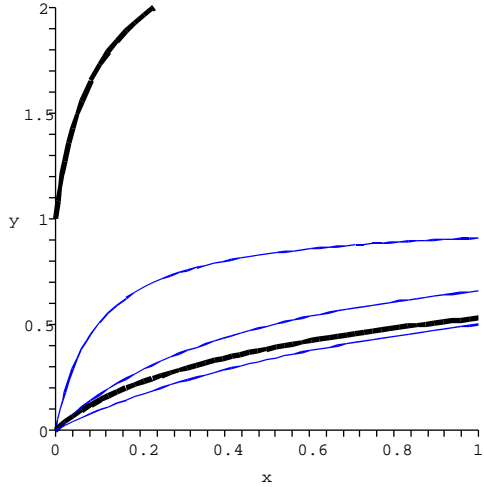}
    \caption{The two thick curves are curves along which solutions
    of \eqref{eq003} have inflection points, for parameter values
    ${\vep=0.6}$ and ${\eta=0.9}$.  The thin curves are the horizontal,
    $\alpha$, and vertical isoclines.} \label{fig005}
\end{center}
\end{figure}

\begin{rem}
We now know that solutions can only have inflection points between
$\cM$ and $\alpha$ and above ${y=1}$.  There are, in fact, curves
along which solutions have zero second derivative.  To find them,
one could, for example, use Maple to solve
${\derivb{f(x,y(x))}{x}=0}$ for $y(x)$, the solutions unfortunately
being rather long and messy. There are three solutions.  One curve
lies below the $x$-axis and is discarded. The other two curves are
in the positive quadrant, one lying between $\cM$ and $\alpha$ (and
which is a lower fence actually), the other starting at $(0,1)$ and
increasing with $x$. See Figure~\ref{fig005}.
\end{rem}

\section{Behaviour of Solutions Near the Origin} \label{sec006}

It was argued in \cite{NguyenFraser} and \cite{Roussel1994}, for
example, that the slow manifold can be written as a Taylor series of
the form ${\cM(x) = \sum_{n=0}^\infty \sigma_n x^n}$ at the origin.
This is a traditional approach but we will show that this approach
is not always valid.  However, in the realm that is usually
considered for the Michaelis-Menten Mechanism, namely ${0 < \vep \ll
1}$, a very high number of terms of this Taylor series is correct.

Intuitively, we know that $\cM$ lies between the horizontal and
vertical isoclines which both have limit zero as ${x\to0^+}$, and
$\cM$ shares the same direction as the slow eigenvector $\vec{v}_+$
at the origin. Hence, it must be that
\begin{subequations} \label{eq016}
\begin{equation} \label{eq016a}
    \sigma_0 = 0
    \mymbox{and}
    \sigma_1 = \sigma.
\end{equation}
By substituting the series into the differential equation, one can
obtain all the coefficients recursively:
\begin{align} \label{eq016b}
    \sigma_n
    &= -\frac{ \sum_{k=2}^{n-1} \sqb{ ( n - k ) \sigma_{n-k} + (1-\eta) (n-k+1) \sigma_{n-k+1} } \sigma_k }{ \vep^{-1} + (1-\eta) (n+1) \sigma_1 - n } \notag\\
    &\quad- \frac{ [ (n-1) \sigma_1 + \vep^{-1} ] \sigma_{n-1} }{ \vep^{-1} + (1-\eta) (n+1) \sigma_1 - n }.
\end{align}
\end{subequations}
Let $y$ be any solution to \eqref{eq003} that lies inside
$\Gamma_0$.  Since no property of the slow manifold was used in
constructing the above series which all other solutions do not
possess, we can equally well write
\begin{equation} \label{eq017}
    y(x) = \sum_{i=0}^\infty \sigma_i x^i.
\end{equation}

Define
\begin{equation} \label{eq018}
    \kappa
    := \frac{ \lambda_- }{ \lambda_+ }
    = \frac{ \vep + 1 + \sqrt{ ( \vep + 1 )^2 - 4 \vep \eta } }{ \vep + 1 - \sqrt{ ( \vep + 1 )^2 - 4 \vep \eta } },
\end{equation}
where we made use of the expression for $\lambda_-$ and $\lambda_+$,
\eqref{eq005}.

Re-arranging \eqref{eq018}, we see that $\eta$ can be written in
terms of $\vep$ and $\kappa$ as
\begin{equation} \label{eq019}
    \eta = \frac{ \kappa ( \vep + 1 )^2 }{ \vep ( \kappa + 1 )^2 }.
\end{equation}
This can tell us when the parameter $\kappa$ takes on certain
values. However, for a given ${\vep>0}$ and
${n\in\Nd{}\backslash\setb{1}}$, there may not be a corresponding
${\eta\in(0,1)}$ that gives ${\kappa=n}$.  It can be shown that
\[
    \kappa = \sfrac{1}{\vep\eta} + \sfrac{2(1-\eta)}{\eta} + \cO(\vep)
    \mymbox{as}
    \vep \to 0
\]
and thus ${\kappa\to\infty}$ as ${\vep \to 0}$.  That is, if $\vep$
is very small, which is the case traditionally considered, $\kappa$
is very large.  Many results that follow will involve $\kappa$ and
so it is a good idea to keep this in mind.

\begin{figure}[t]
\begin{center}
    \includegraphics[width=\mediumfigwid]{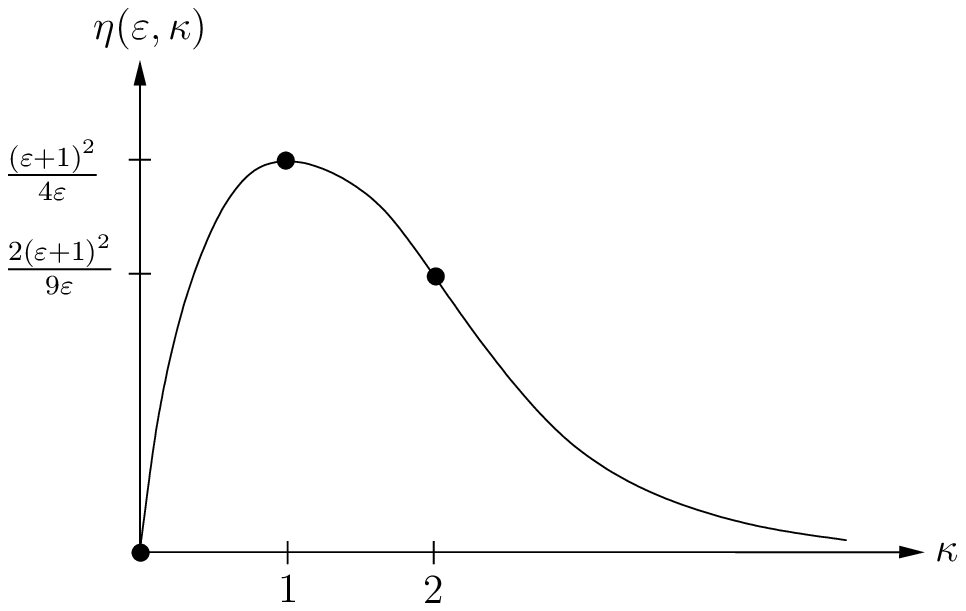}
    \caption{The graph of $\eta(\vep,\kappa)$ for an arbitrary, fixed ${\vep>0}$. The physically relevant values of $\kappa$ are
    ${\kappa>1}$.  Observe that ${\pderiv{\eta}{\kappa}>0}$ for ${0<\kappa<1}$ and ${\pderiv{\eta}{\kappa}<0}$ for ${\kappa>1}$ with a global maximum at ${\kappa=1}$.  Furthermore, there is
    an inflection point at ${\kappa=2}$ and ${\eta \to 0}$ as ${\kappa\to\infty}$. Observe that the maximum value satisfies ${\sfrac{(\vep+1)^2}{4\vep} \geq 1}$ for all ${\vep>0}$ and hence for any permissible value of $\vep$
    there are values of $\kappa$ (in a neighbourhood of ${\kappa=1}$) which give inadmissible values of $\eta$.} \label{fig006}
\end{center}
\end{figure}

Observe that ${\kappa>1}$ and that we can choose values of the
parameters $\vep$ and $\eta$ to achieve any desired value of
$\kappa$ we wish.  The following is easy to prove.

\begin{prop}
Consider the constant $\kappa>1$ and the coefficient $\sigma_2$.
\begin{enumerate}[(a)]
    \item
        There is resonance with the eigenvalues $\setb{\lambda_-,\lambda_+}$ if and only if ${\kappa \in
        \Nd{} \backslash \setb{1}}$.  (See, for example, \cite{Arnold} for a discussion of
        resonance.)
    \item
        For the numbers $\kappa$ and $\sigma_2$, ${\kappa \in (1,2)}$ if and
        only if ${\sigma_2>0}$.  Furthermore, ${\kappa > 2}$ if and only if
        ${\sigma_2<0}$.
\end{enumerate}
\end{prop}

\begin{prop}
Consider the constant $\kappa>1$ and the coefficients
$\seq{\sigma}{i}{0}{\infty}$.
\begin{enumerate}[(a)]
    \item
        If ${\kappa \not\in \Nd{}}$, then all of $\seq{\sigma}{i}{0}{\infty}$ are defined.
    \item
        If ${\kappa \in \Nd{} \backslash \setb{1}}$, then $\seq{\sigma}{i}{0}{\kappa-1}$ are all defined but
        $\sigma_\kappa$ is not defined (and hence all subsequent $\sigma_n$ are not defined).
\end{enumerate}
\end{prop}

\begin{proof}
We know from \eqref{eq016a} that the coefficients $\sigma_0$ and
$\sigma_1$ are always defined.  Consider the expression
\eqref{eq016b}, which gives the recursive descriptions of the
coefficients.  Solving
\[
    \vep^{-1} + ( 1 - \eta ) ( j + 1 ) \sigma_1 - j = 0,
    \quad
    j \in \setb{2,3,\ldots}
\]
gives
\[
    \eta = \frac{ j ( \vep + 1 )^2 }{ \vep ( j + 1 )^2 }.
\]
From \eqref{eq019}, this is true if and only if ${\kappa=j}$. Hence,
if ${\kappa\in\Nd{}\backslash\setb{1}}$, then
$\seq{\sigma}{i}{0}{\kappa-1}$ are all defined but $\sigma_\kappa$
is not defined.  If ${\kappa\not\in\Nd{}\backslash\setb{1}}$, then
all the coefficients are defined.
\end{proof}

A classic method of finding an asymptotic expression for a solution
to a one-dimensional differential equation ${y'=f(x,y)}$ is the
power series method.  Here, one \emph{assumes} a solution of the
form ${y(x)=\sum_{n=0}^\infty a_n x^n}$, substitutes into the
differential equation, and arrives at recursive relationships for
the coefficients which are then solved.  However, this is not always
reliable.

\begin{thm} \label{thm003}
Consider the system of ordinary differential equations
\begin{equation} \label{eq020}
    \dot{ \vec{x} } = \mat{A} \vec{x} + \vec{b}( \vec{x} ),
    \quad
    \vec{x}(0) = \vec{x}_0,
    \quad
    \vec{x} = \stwobyonematrix{x}{y} \in \Rd{2},
\end{equation}
where the matrix $\mat{A}$ is Hurwitz (asymptotically stable), the
vector field $\vec{b}$ is analytic with ${\norm{ \vec{b}(\vec{x})} =
\cO\rb{\norm{\vec{x}}^2}}$ as ${\norm{\vec{x}} \to 0}$, and
$\norm{\vec{x}_0}$ is sufficiently small.  Let the eigenvalues be
$\lambda_+$ and $\lambda_-$, where ${\lambda_-<\lambda_+<0}$, and
define the ratio ${\kappa := \frac{\lambda_-}{\lambda_+} > 1}$.
Suppose that ${\kappa \not\in \Nd{}}$ (i.e. no resonance) and the
eigenvector $\vec{v}_+$ satisfies ${(\vec{v}_+)_1 \ne 0}$. If
$\vec{x}(t)$ is a solution to \eqref{eq020} which approaches the
origin in the slow direction and (for simplicity) is strictly
positive for sufficiently large $t$, then
\[
    y(t) = \sum_{n=1}^{\floor{\kappa}} \sigma_n x(t)^n + C x(t)^\kappa + o\rb{ x(t)^\kappa }
    \mymbox{as}
    t \to \infty
\]
for some constants $\seq{\sigma}{n}{1}{\floor{\kappa}}$ (which are
independent of initial condition) and $C$ (which depends on the
initial condition).
\end{thm}

\begin{proof}
In order to derive the necessary asymptotic expansion for $y(t)$ in
terms of $x(t)$, we make use of the linearized problem
\begin{equation} \label{eq021}
    \dot{\vec{z}} = \mat{A} \vec{z},
    \quad
    \vec{z}(0) = \vec{z}_0.
\end{equation}
To avoid the trivial solutions, which have nothing to offer us, we
will assume that ${\vec{x}_0,\vec{z}_0 \ne \vec{0}}$. Let
$\vec{x}(t)$ and $\vec{z}(t)$ be, respectively, the unique solutions
to \eqref{eq020} and \eqref{eq021}, both of which tend to the origin
as time tends to infinity. We will not consider the initial
conditions $\vec{x}_0$ and $\vec{z}_0$ to be independent so that the
solutions $\vec{x}(t)$ and $\vec{z}(t)$ can be related. Furthermore,
we need both $\norm{\vec{x}_0}$ and $\norm{\vec{z}_0}$ to be small.

The solution to the linear problem $\vec{z}(t)$ can be written in
the explicit form
\[
    \vec{z}(t) = c_+ \exp{\lambda_+ t} \vec{v}_+ + c_- \exp{\lambda_- t} \vec{v}_-,
\]
where ${c_+>0}$ (since we assumed that solutions approach the origin
from the right in the slow direction).

We know that there is no resonance with the eigenvalues.  Moreover,
the eigenvalues are in the Poincar\'{e} domain. Applying
Poincar\'{e}'s Theorem (see, for example, \cite{Arnold} p.190),
there is a quadratic vector field $\vec{q}$ such that ${\vec{x} =
\vec{z} + \vec{q}(\vec{z})}$. Hence, we can write (not uniquely if
${\kappa \in \Qd{}}$)
\begin{subequations} \label{eq022}
\begin{align}
    x(t)
    &\sim \sum_{(m,n) \in S} a_{mn} \exp{ ( m \lambda_- + n \lambda_+ ) t }
    = \sum_{(m,n) \in S} a_{mn} \exp{ ( m \kappa + n ) \lambda_+ t } \mymbox{as} t \to \infty \label{eq022a}\\
    y(t)
    &\sim \sum_{(m,n) \in S} b_{mn} \exp{ ( m \lambda_- + n \lambda_+ ) t }
    = \sum_{(m,n) \in S} b_{mn} \exp{ ( m \kappa + n ) \lambda_+ t } \mymbox{as} t \to \infty, \label{eq022b}
\end{align}
\end{subequations}
where
\[
    S := \setb{ (m,n) \st m,n \in \Zd{}, ~ m,n \geq 0, ~ m + n \geq 1 }.
\]

Let ${\ell := \floor{\kappa}}$.  Then, the first ${\ell+1}$ most
dominant terms in \eqref{eq022} are, in order of decreasing
dominance,
\[
    \exp{\lambda_+ t},
    \exp{2 \lambda_+ t},
    \ldots,
    \exp{\ell \lambda_+ t},
    \exp{\lambda_- t}.
\]
To see why this is the case, we make two observations. First, the
fact that the listed exponentials are in decreasing order of
dominance is obvious except maybe for the last two.  Since ${\kappa
= \frac{\lambda_-}{\lambda_+} > \ell}$, we have ${\lambda_- < \ell
\lambda_+ < 0}$.  Second, there cannot be any other exponentials of
the form $\exp{ ( m \lambda_+ + n \lambda_- ) t}$ in between those
listed.

For our purposes, we need only the first ${\ell+1}$ terms of
\eqref{eq022} and hence we write
\begin{subequations} \label{eq023}
\begin{align}
    x(t) &= \sum_{m=1}^\ell a_m \exp{m \lambda_+ t} + a_{\ell+1} \exp{\lambda_- t } + o\rb{ \exp{\lambda_- t } } \mymbox{as} t \to \infty \label{eq023a} \\
    \quad
    y(t) &= \sum_{m=1}^\ell b_m \exp{m \lambda_+ t} + b_{\ell+1} \exp{\lambda_- t } + o\rb{ \exp{\lambda_- t } } \mymbox{as} t \to \infty, \label{eq023b}
\end{align}
\end{subequations}
where ${a_1 \ne 0}$. The coefficients can be related using the
differential equation.  To write $y(t)$ in terms of $x(t)$, we will
successively eliminate the exponentials. Manipulating \eqref{eq023a}
and \eqref{eq023b},
\begin{equation} \label{eq024}
    y(t) - \sigma_1 x(t)
    = \sum_{m=2}^\ell b_m^{(2)} \exp{m \lambda_+ t} + b_{\ell+1}^{(2)} \exp{\lambda_- t } + o\rb{ \exp{\lambda_- t } },
\end{equation}
where ${\sigma_1 := \sfrac{ b_1 }{ a_1 }}$ and ${b_m^{(2)} := b_m -
\sigma_1 a_m}$. To go further, observe that we can write powers of
$x(t)$ as
\[
    x(t)^n = a_1^n \exp{n \lambda_+ t} + \sum_{m=n+1}^\ell c_{mn} \exp{m \lambda_+ t} + o\rb{ \exp{\lambda_- t } }
    \mymbox{as}
    t \to \infty,
\]
where ${n \in \setb{ 2, \ldots, \ell }}$. Solving for the most
dominant exponential,
\begin{equation} \label{eq025}
    \exp{n \lambda_+ t} = \sfrac{ 1 }{ a_1^n } x(t)^n - \sum_{m=n+1}^\ell \rb{ \sfrac{ c_{mn} }{ a_1^n } } \exp{m \lambda_+ t} + o\rb{ \exp{\lambda_- t } }
    \mymbox{as}
    t \to \infty.
\end{equation}
With \eqref{eq025}, we can successively eliminate the exponentials
of \eqref{eq024}---each time introducing other exponential terms but
none of order already eliminated---until we are left with an
expression of the form
\begin{equation} \label{eq026}
    y(t) - \sum_{m=1}^\ell \sigma_m x(t)^m = b_{\ell+1}^{(\ell+1)} \exp{\lambda_- t} + o\rb{ \exp{\lambda_- t } }
    \mymbox{as}
    t \to \infty.
\end{equation}
Since ${x(t)^\kappa = a_1^\kappa \exp{\lambda_- t} + o\rb{
\exp{\lambda_- t} }}$ and hence
\[
    \exp{\lambda_- t} = \sfrac{1}{a_1^\kappa} x(t)^\kappa + o\rb{ \exp{\lambda_- t} }
    \mymbox{as}
    t \to \infty,
\]
we can write \eqref{eq026} as
\[
    y(t) - \sum_{m=1}^\ell \sigma_m x(t)^m = C x(t)^\kappa + o\rb{
    \exp{\lambda_- t} }
    \mymbox{as}
    t \to \infty,
\]
where
\[
    C := \frac{ b_{\ell+1} - \sigma_1 a_{\ell+1} }{ a_1^\kappa }.
\]
The desired conclusion follows.
\end{proof}

\begin{rem}
The coefficients $\seq{\sigma}{n}{1}{\ell}$ are calculated using the
power series method.  That is, one assumes that the solution to the
one-dimensional version of \eqref{eq020} is ${y(x) =
\sum_{n=0}^\infty \sigma_n x^n}$ (where ${\sigma_0=0}$ by
necessity).  The purpose of the theorem is to tell us how many of
the resulting terms apply to all solutions.
\end{rem}

\begin{cor} \label{cor001}
There exists a solution to \eqref{eq020} such that
\[
    y(t) \sim \sum_{n=1}^\infty \sigma_n x(t)^n
    \mymbox{as}
    t \to \infty
\]
for some constants $\seq{\sigma}{n}{1}{\infty}$.
\end{cor}

\begin{proof}
Choose the initial condition $\vec{z}_0$ so that it is parallel to
the slow eigenvector $\vec{v}_+$.  Then,
\[
    \vec{z}(t) = c_+ \exp{\lambda_+ t} \vec{v}_+
\]
and hence we can write
\[
    x(t) \sim \sum_{m=1}^\infty a_m \exp{m \lambda_+ t}
    \mymbox{and}
    y(t) \sim \sum_{m=1}^\infty b_m \exp{m \lambda_+ t} \mymbox{as} t \to \infty.
\]
Any positive integer power of $x(t)$ will be a series of the same
form as $x(t)$ and $y(t)$.  Successively eliminating exponents, just
like in the proof of the theorem, gives us our desired conclusion.
\end{proof}

Now, we apply this general result to the Michaelis-Menten Mechanism.

\begin{lem} \label{lem002}
Let $y$ be a solution to \eqref{eq003} lying inside $\Gamma_0$ and
let the ratio of the eigenvalues be ${\kappa >1}$.  Suppose that
${\kappa \not\in \Nd{}}$ (i.e. there is no resonance).  Then,
\[
    y(x) = \sum_{n=1}^{\floor{\kappa}} \sigma_n x^n + C x^\kappa + o\rb{ x^\kappa }
    \mymbox{as}
    x \to 0^+,
\]
where $\seq{\sigma}{n}{1}{\floor{\kappa}}$ are as in \eqref{eq017}
and $C$ is some constant that distinguishes the solution $y(x)$ from
other such solutions.
\end{lem}

\begin{proof}
This proof is a simple application of Theorem~\ref{thm003}.  The
original, time-dependent differential equation \eqref{eq001} can be
written
\begin{equation} \label{eq027}
    \dot{\vec{x}} = \mat{A} \vec{x} + \vec{b}(\vec{x}),
    \quad
    \vec{x}(0) = \vec{x}_0
\end{equation}
where the analytic vector field $\vec{b}$ is given by
\[
    \vec{b}(\vec{x}) := x y \stwobyonematrix{1}{-\vep^{-1}}.
\]
Observe that
\[
    \norm{ \vec{b(\vec{x})} } = \cO\rb{ \norm{\vec{x}}^2 }
    \mymbox{as}
    \norm{\vec{x}} \to 0.
\]
This is fairly obvious but can be shown directly:
\[
    0
    \leq
    ( x - y )^2
    = x^2 + y^2 - 2 x y
    = \norm{\vec{x}}^2 - \sfrac{2}{\sqrt{1+\vep^{-2}}} \norm{\vec{b}(\vec{x})}
\]
and hence
\[
    \norm{\vec{b}(\vec{x})} \leq \sfrac{\sqrt{1+\vep^{-2}}}{2} \norm{\vec{x}}^2.
\]
The linearized problem is
\begin{equation} \label{eq028}
    \dot{\vec{z}} = \mat{A} \vec{z},
    \quad
    \vec{z}(0) = \vec{z}_0,
\end{equation}
where the matrix
${\mat{A}=\stwobytwomatrix{-1}{1-\eta}{\vep^{-1}}{-\vep^{-1}}}$ was
given in \eqref{eq004}. Again, we will assume that
${\vec{x}_0,\vec{z}_0 \ne \vec{0}}$ and, in particular, lie in the
positive quadrant. Let $\vec{x}(t)$ and $\vec{z}(t)$ be,
respectively, the unique solutions to \eqref{eq027} and
\eqref{eq028}.  The initial conditions $\vec{x}_0$ and $\vec{z}_0$
are not independent so that the solutions $\vec{x}(t)$ and
$\vec{z}(t)$ can be related. Furthermore, we need both
$\norm{\vec{x}_0}$ and $\norm{\vec{z}_0}$ to be small. The solution
to the linear problem $\vec{z}(t)$, as we have seen earlier, can be
written explicitly as
\[
    \vec{z}(t) = c_+ \exp{\lambda_+ t} \vec{v}_+ + c_- \exp{\lambda_- t} \vec{v}_-,
    \quad
    c_\pm := \frac{ \hvec{v}_\pm^T \vec{z}_0 }{ \hvec{v}_\pm^T \vec{v}_\pm }
\]
with ${c_+>0}$ and the sign of $c_-$ depending on which side of
$\vec{v}_+$ the initial point $\vec{z}_0$ lies.

Finally applying Theorem~\ref{thm003}, after dropping the time
dependence we can say that
\[
    y(x) = \sum_{n=1}^\ell \hsigma_n x^n + C x^\kappa + o\rb{x^\kappa}
    \mymbox{as} x \to 0^+
\]
for some constants $\seq{\hsigma}{n}{1}{\ell}$ and $C$, where ${\ell
:= \floor{\kappa}}$. By uniqueness, we have ${\hsigma_n=\sigma_n}$
for all ${n\in\setb{1,\ldots,\ell}}$.
\end{proof}

\begin{rem}
For ${\kappa\in(1,2)}$, we can manipulate the given constants to get
\[
    C
    = \frac{ c_- ( \lambda_- - \lambda_+ ) }{ c_+^\kappa ( 1 - \eta )^\kappa }.
\]
\end{rem}

\section{All or Most Solutions Must Enter the Antifunnel} \label{sec007}

We now investigate conditions under which solutions enter
$\Gamma_1$.

\begin{thm} \label{thm004}
Let $\vec{x}(t)$ be a solution to \eqref{eq001} and suppose there is
no resonance, i.e. ${\kappa \not\in \Nd{}}$.
\begin{enumerate}[(a)]
    \item\label{thm004a}
        If ${\kappa>2}$, then there exists a $t^*>0$ such that
        \[
            (x(t),y(t)) \in \Gamma_1 ~ \forall ~ t \geq t^*.
        \]
    \item\label{thm004b}
        If ${\kappa < 2}$, then there exist solutions $\vec{x}(t)$ which do not enter $\Gamma_1$.
        Moreover, solutions that do not enter $\Gamma_1$ must
        enter $\Gamma_0$ through the vertical isocline $V$ to the
        left of the line ${y = \sigma x}$.
\end{enumerate}
\end{thm}

\begin{figure}[t]
\begin{center}
    \includegraphics[width=\mediumfigwid]{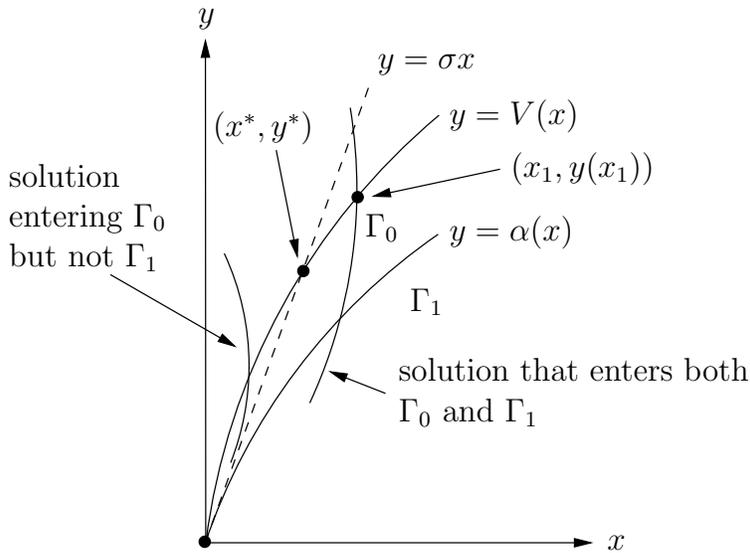}
    \caption{Sketch demonstrating Theorem~\ref{thm004}(b).} \label{fig007}
\end{center}
\end{figure}

\begin{proof}
We begin by noting that if a solution $\vec{x}(t)$ enters
$\Gamma_1$, it forever remains in $\Gamma_1$.  This is because
${\vec{g}\dotprod\normal<0}$ along $\alpha$ and $H$, where $\normal$
is the unit normal vector.  Let $y$ be the corresponding
one-dimensional solution.

\begin{enumerate}[(a)]
    \item
        Applying Lemma~\ref{lem002}, we can write
        \[
            y(x) = \sigma x + \sigma_2 x^2 + o(x^2)
            \mymbox{as}
            x \to 0^+.
        \]
        Furthermore, since ${\kappa>2}$ we have ${\sigma_2<0}$. Thus,
        \[
            \lim_{x \to 0^+} y''(x) = 2 \sigma_2 < 0.
        \]
        Since solutions are concave down only when they lie below
        the isocline $\alpha$, $\vec{x}(t)$ eventually enters the $\Gamma_1$ antifunnel.
    \item
        From Lemma~\ref{lem002} we have
        \[
            y(x) = \sigma x + C x^\kappa + o\rb{x^\kappa}
            \mymbox{as}
            x \to 0^+
        \]
        for some constant $C$.  It follows that there are some solutions to \eqref{eq003} that are concave
        up at the origin---the ones for which ${C>0}$---and curve away from $\alpha$ and exit $\Gamma_0$
        through the vertical isocline. Moreover, since ${\sigma_2>0}$, by virtue of Corollary~\ref{cor001}
        it follows that there is a solution with a Taylor series at the origin that does not enter $\Gamma_1$
        from above.  See Figure~\ref{fig007}.

        We know already that $\vec{x}(t)$ eventually enters $\Gamma_0$. If $\vec{x}(t)$ enters
        $\Gamma_0$ through the horizontal isocline, it also enters
        $\Gamma_1$.  Denote by $(x^*,y^*)$ the point of intersection of the line ${y = \sigma x}$ and the vertical isocline ${y=V(x)}$
        (that is, ${y^* = \sigma x^* = V(x^*)}$).
        Assume that $\vec{x}(t)$ enters $\Gamma_0$
        through the vertical isocline to the right of $(x^*,y^*)$.  We claim that $\vec{x}(t)$ also
        enters $\Gamma_1$. Suppose that $y$ intersects the vertical isocline at ${x=x_1}$.  Observe that
        ${x_1 \geq x^*}$, by assumption, and ${y(x_1) \leq \sigma
        x_1}$. By the Mean Value Theorem, there is a ${x_0\in(0,x_1)}$ such that
        \[
            0
            \leq y'(x_0)
            = \frac{ y(x_1) - 0 }{ x_1 - 0 }
            = \frac{y(x_1)}{x_1}
            \leq \frac{ \sigma x_1 }{ x_1 }
            = \sigma.
        \]
        By Virtue of the isocline structure, ${H(x_0) \leq y(x_0) \leq \alpha(x_0)}$ and therefore $\vec{x}(t)$ eventually enters
        $\Gamma_1$.
\end{enumerate}
\end{proof}

\section{Properties of the Slow Manifold} \label{sec008}

Finally, we present some important properties of the slow manifold.

\begin{prop}
The slow manifold ${y=\cM(x)}$ is concave down for all $x>0$.
\end{prop}

\begin{proof}
Construct a sequence of functions $\seq{y}{n}{N}{\infty}$ as
follows.  Fix ${x_0>0}$ and let $y_n$ be the solution to
\eqref{eq003} such that
\[
    y_n(x_0) = \cM(x_0) - \sfrac{1}{n}.
\]
The number $N$ is taken large enough so that ${y_N(x_0) \geq 0}$.
Thus, by Proposition~\ref{prop002}, ${y_n''(x_0) < 0}$ for all $n$.
Now,
\begin{equation} \label{eq029}
    y_n''(x_0) = p(x_0,y_n(x_0)) h(x_0,y_n(x_0))
\end{equation}
and
\begin{equation} \label{eq030}
    \cM''(x_0) = p(x_0,\cM(x_0)) h(x_0,\cM(x_0)).
\end{equation}
By construction,
\[
    \lim_{n\to\infty} y_n(x_0) = \cM(x_0).
\]
Letting ${n\to\infty}$ in \eqref{eq029} and applying \eqref{eq030},
continuity tells us
\[
    \lim_{n\to\infty} y_n''(x_0) = \cM''(x_0).
\]
Since ${y_n''(x_0)<0}$ for all $n$, ${\cM''(x_0) \leq 0}$.  Since
$x_0$ was arbitrary, ${\cM''(x) \leq 0}$ for all ${x>0}$.

We now establish a strict inequality.  Suppose that ${\cM''(x^*)=0}$
for ${x^*>0}$.  If $h$ is as in \eqref{eq013} with respect to the
solution $\cM$, we have $h'(x^*)>0$.  This contradicts the fact that
${h(x) \leq 0}$ for all $x>0$.
\end{proof}

\begin{prop}
The slow manifold ${y=\cM(x)}$ satisfies, for all $x>0$,
\[
    0 < H(x) < \cM(x) < \alpha(x) < 1.
\]
Furthermore,
\[
    \lim_{x \to 0^+} \cM(x) = 0
    \mymbox{and}
    \lim_{x\to\infty} \cM(x) = 1.
\]
\end{prop}

\begin{proof}
The first result follows from Theorem~\ref{thm002} and the
definition of $\alpha$, \eqref{eq011}.  To establish the other two
results, note that the functions $H$, $\alpha$, and $\cM$ are all
continuous for $x>0$. Since
\[
    \lim_{x \to 0^+} H(x) = \lim_{x \to 0^+} \alpha(x) = 0
    \mymbox{and}
    \lim_{x \to \infty} H(x) = \lim_{x \to \infty} \alpha(x) = 1,
\]
the results follow from the Squeeze Theorem.
\end{proof}

\begin{prop}
The slope of the slow manifold $y=\cM(x)$ satisfies, for ${x>0}$,
\[
    0 < \cM'(x) < \sigma.
\]
Furthermore,
\[
    \lim_{x \to 0^+} \cM'(x) = \sigma
    \mymbox{and}
    \lim_{x\to\infty} \cM'(x) = 0.
\]
\end{prop}

\begin{proof}
The first result follows from the fact that the slow manifold lies
within $\Gamma_1$, which consists of nested isoclines of slopes
varying from $0$ to $\sigma$.

To prove the second result, observe that the direction of the slow
manifold at the origin must correspond to the slow eigenvector at
the origin, which has slope $\sigma$.

To prove the third result, we note that, since $\cM$ is strictly
increasing and concave down, there is a ${c \in [0,\sigma)}$ such
that
\[
    \lim_{x \to \infty} \cM'(x) = c.
\]
Suppose, on the contrary, that $c>0$.  By the Fundamental Theorem of
Calculus,
\[
    \cM(x)
    = \integral{0}{x}{\cM'(u)}{u}
    > \integral{0}{x}{c}{u}
    = c x.
\]
However, for sufficiently large $x$, ${cx > \alpha(x)}$, a
contradiction.
\end{proof}

The slow manifold has also been approximated, for large $x$, in the
asymptotic series\cite{Roussel1994}
\begin{equation} \label{eq031}
    \cM(x) \sim \rho_0 + \rho_1 x^{-1} + \rho_2 x^{-2} + \cdots
    \mymbox{as}
    x \to \infty.
\end{equation}
The coefficients can be obtained by substituting the series into the
differential equation and are given recursively by
\begin{align} \label{eq032}
    \rho_0 &= 1, \quad \rho_1 = -1, \quad \rho_2 = 1, \notag \\
    \rho_{n} &= -\rho_{n-1} + \vep \sum_{i=1}^{n-2} i \rho_i [ \rho_{n-i-1} + ( 1 - \eta ) \rho_{n-i-2} ]
    \mymbox{for} n > 2.
\end{align}
Observe that all the coefficients are polynomials in $\vep$ and
$\eta$.  As we will establish now, the series \eqref{eq031} is fully
correct.

\begin{prop}
For large $x$, the slow manifold satisfies
\[
    \cM(x) \sim \sum_{n=0}^\infty \rho_n x^{-n}
    \mymbox{as}
    x \to \infty
\]
where $\seq{\rho}{i}{0}{\infty}$ are as in \eqref{eq032}. For small
$x$, in the case of no resonance (i.e. ${\kappa \not\in \Nd{}}$),
the slow manifold satisfies
\[
    \cM(x) = \sum_{i=1}^{\floor{\kappa}} \sigma_i x^i + C x^\kappa + o\rb{x^\kappa}
    \mymbox{as}
    x \to 0^+
\]
for some constant $C$.
\end{prop}

\begin{proof}
The second conclusion follows from Lemma~\ref{lem002}.  To prove the
first conclusion, observe that for any ${c>0}$ there exists a
${x^*>0}$ such that
\[
    H(x) < \cM(x) < F(x,c)
\]
for all ${x>x^*}$. This is because $\cM$ is concave down and
${\lim_{x\to\infty} \cM'(x) = 0}$. Hence,
\[
    1 - x^{-1} + \cO\rb{x^{-2}} < \cM(x) < 1 - K(c) x^{-1} + \cO\rb{x^{-2}}
    \mymbox{as}
    x \to \infty.
\]
Since ${c>0}$ was arbitrary and ${K(c) \to 1}$ as ${c \to 0}$, it
follows that
\[
    \lim_{x\to\infty} \frac{ \cM(x) - 1 }{ -x^{-1} } = 1.
\]
Hence,
\[
    \cM(x) = 1 - x^{-1} + o\rb{x^{-1}}
    \mymbox{as}
    x \to \infty.
\]
Unfortunately, this is as much information that we can extract using
the isoclines.  To obtain the remaining terms of the asymptotic
series, we will use the Centre Manifold Theorem (see, for example,
\cite{Carr}).

Under the change of variables
\begin{equation} \label{eq033}
    X := x^{-1},
    \quad
    Y := y - \rb{ 1 - x^{-1} },
\end{equation}
we arrive at the system
\begin{align} \label{eq034}
    \dot{X} &= -X^2 g_1\rb{ X^{-1}, 1 - X + Y } \notag\\
    \dot{Y} &= -X^2 g_1\rb{ X^{-1}, 1 - X + Y } + g_2\rb{ X^{-1}, 1 - X + Y },
\end{align}
where $g_1$ and $g_2$ are as in \eqref{eq002}.  The system
\eqref{eq034} is not polynomial but there is no harm, because the
resulting one-dimensional differential equation will be the same, in
considering the system
\begin{align} \label{eq035}
    \dot{X} &= -X^3 g_1\rb{ X^{-1}, 1 - X + Y } \notag\\
    \dot{Y} &= -X^3 g_1\rb{ X^{-1}, 1 - X + Y } + X g_2\rb{ X^{-1}, 1 - X + Y },
\end{align}
which is polynomial.  Expanding, we get the expressions
\begin{align*}
    X^3 g_1\rb{ X^{-1}, 1 - X + Y } &= -X^2 \sqb{ \eta X - Y + \rb{ 1 - \eta } X \rb{ X - Y  } } \\
    X g_2\rb{ X^{-1}, 1 - X + Y } &= \vep^{-1} \rb{ X^2 - X Y - Y }.
\end{align*}
The system \eqref{eq035}, as we see, is a bit more messy than the
original system \eqref{eq003}. The eigenvalues of the matrix for the
linear part of the new system \eqref{eq035}, which is diagonal by
construction, are $0$ and $-\vep^{-1}$. We know from centre manifold
theory that there is a centre manifold which, we claim, must be the
slow manifold.

Observe that the $Y$-axis is invariant. Moreover, the fixed point
${(X,Y)=(0,0)}$ is a saddle node (or a degenerate saddle).  The
physically relevant portion of the phase portrait, namely ${X \geq
0}$ and ${Y \geq -1}$, consists of two hyperbolic sectors, one with
the positive $Y$-axis and the centre manifold as boundaries and the
other with the negative $Y$-axis and the centre manifold as
boundaries. See Figure~\ref{fig008}. This can be shown using
techniques in \S9.21 of \cite{AndronovLeontovichGordonMaier} (in
particular Theorem~65).  This also is a consequence of the phase
portrait of the original system \eqref{eq001}. It follows that the
centre manifold is indeed the slow manifold.

\begin{figure}[t]
\begin{center}
    \includegraphics[width=\largefigwid]{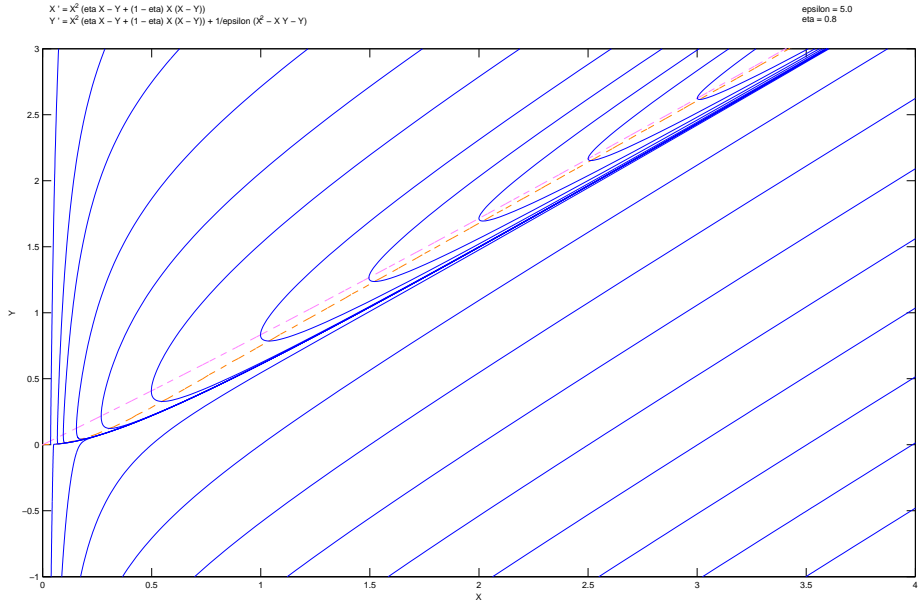}
    \caption{A phase portrait for \eqref{eq035} for ${\vep=5.0}$
    and ${\eta=0.8}$.} \label{fig008}
\end{center}
\end{figure}

By the Centre Manifold Theorem, the slow manifold (in the new
coordinates) can be written
\[
    \cM(X) \sim \sum_{i=2}^\infty \hrho_i X^i
    \mymbox{as}
    X \to 0^+,
\]
for some coefficients $\seq{\hrho}{i}{2}{\infty}$.  The first two
terms are given by ${\hrho_2=1}$ and ${\hrho_3 = \vep \eta - 1}$.
Reverting back to the original coordinates,
\[
    \cM(x) \sim 1 - x^{-1} + \sum_{i=2}^\infty \hrho_i x^{-i}
    \mymbox{as}
    x \to \infty.
\]
Observing that the coefficients in \eqref{eq032} are generated
uniquely from the differential equation, the conclusion follows.
\end{proof}

\needspace{2.0cm}
\begin{rems}~
\begin{enumerate}[(i)]
    \item
        We use the ad hoc transformation \eqref{eq033} because it is inspired by the
        series for $\cM$ which we wish to obtain and it also results in a
        system which is in the canonical form of the Centre Manifold
        Theorem. Others, for example,
        \cite{DavisSkodje,GearKaperKevrekidisZagaris}, have used
        Poincar\'{e} compactification to study the behaviour of $\cM$ at
        infinity and found that the fixed point is a degenerate saddle.
    \item
        The Centre Manifold Theorem can be applied at the origin as well.
        See, for example, \cite{Carr} pages 8--10.  However, this result
        gives a smooth solution for small $\vep$ only.  This is because
        in order to apply the Centre Manifold Theorem, the differential equation ${\dot{\vep}=0}$
        is appended to the system \eqref{eq002} which gives the zero eigenvalue. A centre
        manifold exists in a neighbourhood of ${(x,y,\vep)=(0,0,0)}$.
\end{enumerate}
\end{rems}

\begin{prop}
In the case of no resonance, i.e. ${\kappa \not\in \Nd{}}$, the
second derivative of the slow manifold satisfies
\[
    \lim_{x \to 0^+} \cM''(x)
    =
    \begin{cases}
        2 \sigma_2, \mymbox{if} \kappa > 2 \\
        -\infty, \mymbox{if} \kappa < 2
    \end{cases}.
\]
\end{prop}

\begin{proof}
The proof involves an easy application of Lemma~\ref{lem002} and the
fact that the slow manifold is concave down at the origin.
\end{proof}

\begin{rem}
The quasi-steady-state approximation has traditionally been used to
approximate the long-term behaviour (in time) of solutions to
\eqref{eq001}.  Is this justified?  Recall that
\[
    H(x) = \frac{x}{1+x}
    \mymbox{and}
    \alpha(x) = \frac{x}{\sigma^{-1}+x}.
\]
It follows that the QSSA is good when ${\sigma \approx 1}$.  Recall
also that ${\sigma = 1 + \cO\rb{\vep}}$ as ${\vep \to 0}$. Hence,
the QSSA is a good approximation when $\vep$ is small.  However, the
function $\alpha(x)$ has slope ${\alpha'(0)=\sigma}$ at the origin
and so is a good approximation for solutions near the origin for any
$\vep$. Furthermore, since ${\kappa = ( \vep \eta )^{-1} + \cO(1)}$
as ${\vep \to 0}$, a large number of Taylor coefficients are correct
in the asymptotic expansion at the origin for the slow manifold if
$\vep$ is small.
\end{rem}

\section{Open Questions} \label{sec009}

In the analysis of the behaviour at the origin, we assume
non-resonance of the eigenvalues of the linearization at the origin.
The resonance cases still need to be investigated.

S. Fraser and M.
Roussel\cite{Fraser1988,NguyenFraser,Roussel1994,Roussel1997,RousselFraser2001}
have introduced and investigated an iteration scheme to approximate
the slow manifold. Specifically, the iterates are defined by
${y_{n+1}:=F(x,y_n')}$, where $F$ is as in \eqref{eq008}. A
definitive proof of convergence of the scheme which is valid for all
values of the parameters $\vep$ and $\eta$ has not yet been given.
The scheme can diverge for certain values of the parameters and
certain choices of initial iterate and converge for others.
Convergence has been examined, for example, in \cite{Kaper}. In this
particular paper, the Fraser iterates and perturbation series for
$\cM$ in $\vep$ were compared. Specifically, if
$\seq{y}{n}{0}{\infty}$ are the Fraser iterates with initial iterate
${y_0:=H}$ and
\[
    \cM(x) = \sum_{m=0}^\infty \cM_m(x) \vep^m
\]
is the perturbation series for $\cM$, then for each $n$
\[
    y_n(x) = \sum_{m=0}^n \cM_m(x) \vep^m + \cO\rb{\vep^{n+1}}
    \mymbox{as}
    \vep \to 0.
\]

\section*{Acknowledgement}

\addcontentsline{toc}{section}{Acknowledgements}

We wish to thank the referee who provided detailed comments and suggestions which led many improvements in exposition and clarity.

\end{document}